\documentclass[11pt]{article}
\usepackage{amsmath,amssymb}
\usepackage{amsfonts}
\usepackage{eucal}
\title{\leavevmode\vadjust{\vskip -5mm}
\textbf{The Group of Parallel Transports \\
in the Riemannian Space}}
\author{\sc Serhiy E. SAMOKHVALOV\thanks{{\bf e}-{\it mail}: samokhval@dstu.dp.ua}\\
 \tabaddress{Department of Applied Mathematics \\
State Technical University, Dniprodzerzhinsk, Ukraine}}

\date{October 10, 2005}
\def\tabaddress#1{{\small\it\begin{tabular}[t]{c}#1
\\[1.2ex]\end{tabular}}}
\def\<#1>{\langle#1\rangle}
\newtheorem{proposition}{Proposition}
\newtheorem{definition}{Definition}

\def\beq{\begin{equation}}
\def\eeq{\end{equation}}
\def\bea{\begin{eqnarray}}
\def\eea{\end{eqnarray}}
\def\beann{\begin{eqnarray*}}
\def\eeann{\end{eqnarray*}}
\def\ben{\begin{enumerate}}
\def\een{\end{enumerate}}
\def\qed{\ifvmode\removelastskip\fi
{\unskip\nobreak\hfil\penalty50\hbox{}\nobreak\hfil \hbox{\vrule
height1.2ex width1.2ex}\parfillskip=0pt \finalhyphendemerits=0
\par\smallskip}}
\def\texthook{\vrule height 0pt depth 0.4pt width 3.5pt
          \vrule height 5pt depth 0.4pt \kern 3pt}
\def\scripthook{\vrule height 0pt depth 0.2pt width 1.5pt
                \vrule height 3pt depth 0.2pt width 0.2pt \kern 1pt}

\begin{document}

\maketitle

\thispagestyle{empty}

\begin{abstract}
\noindent We show that there is an infinite group of special
automorphisms of the deformed group of diffeomorphisms, which
describes parallel transports in Riemannian spaces of any variable
curvature. Generators of translations of such group contain
covariant derivatives, and structure functions - the curvature
tensor.
\end{abstract}


\\
%{\sl PACS number (2003)}:02.40.Yy, 45.20.Jj

%%%%%%%%%%%%%%%%%%%%%%%%%%%%%%%%%%%%%%%%%%%%%%%%%%%%%%%%%%%%%%%%
\clearpage
%%%%%%%%%%%%%%%%%%%%%%%%%%%%%%%%%%%%%%%%%%%%%%%%%%%%%%%%%%%%%%%%
%\tableofcontents
%%%%%%%%%%%%%%%%%%%%%%%%%%%%%%%%%%%%%%%%%%%%%%%%%%%%%%%%%%%%%%%%
%\clearpage
%%%%%%%%%%%%%%%%%%%%%%%%%%%%%%%%%%%%%%%%%%%%%%%%%%%%%%%%%%%%%%%%

%%%%%%%%%%%%%%%%%%%%%%%%%%%%%%%%%%%%%%%%%%%%%%%%%%%%%%%%%%%%%%%%
\section{Introduction}
%%%%%%%%%%%%%%%%%%%%%%%%%%%%%%%%%%%%%%%%%%%%%%%%%%%%%%%%%%%%%%%%

On the way of realization of the Klein's Erlangen Program \cite{1}
for the geometrical structure of a Riemannian space with arbitrary
variable curvature in the paper \cite{2} there is shown, that
Riemannian structure in a manifold $M$ may be naturally set by the
deformed group of its diffeomorphisms $T^{gH}_M = Diff\ M$.
Information about geometrical structure is contained in the
multiplication law of the group $T^{gH}_M$, which sets the rule of
parallel transports of vectors in the tangent bundle $TM$ over
$M$.

In this article we shall show that \emph{parallel transports of
vectors in Riemannian spaces of any variable curvature is
described by an infinite deformed group $DT$ of special
automorphisms of the deformed group of diffeomorphisms $T^{gH}_M$,
which sets this geometrical structure on manifold $M$.} Here we
construct this group and describe some of its properties.
Specifically we show that the covariant derivatives are among
generators and the Riemann curvature tensor is among structure
functions of the group $DT$.

%%%%%%%%%%%%%%%%%%%%%%%%%%%%%%%%%%%%%%%%%%%%%%%%%%%%%%%%%%%%%%%%
\section{The Group of Deformations of the
Generalized Gauge Groups} \protect\label{GDGGG}
%%%%%%%%%%%%%%%%%%%%%%%%%%%%%%%%%%%%%%%%%%%%%%%%%%%%%%%%%%%%%%%%

Let's consider a local Lie group $G_M$ with parameters
$\tilde{g}^\alpha$ (indices $\alpha$, $\beta$, $\gamma$, $\delta$)
and the multiplication law
$(\tilde{g}\cdot\tilde{g}')^{\alpha}=\tilde{\varphi}^{\alpha}
(\tilde{g},\tilde{g}')$, which acts (perhaps inefficiently) on a
coordinate chart $U$ of the manifold $M$ with coordinates $x^\mu$
(indices $\mu$, $\nu$, $\pi$, $\rho$, $\sigma$) according to the
formula $x'^{\mu}=\tilde{f}^{\mu}(x, \tilde{g})$. The local
infinite Lie group $G^g_M$ is parameterized by smooth functions
$\tilde{g}^{\alpha}(x)$ which satisfy the condition
\begin{equation}\label{eq1}
\det\{d_{\nu}\tilde{f}^{\mu}(x, \tilde{g}(x))\}\neq0\quad\forall
x\in U,
\end{equation}
where $d_{\nu}:=d/dx^{\nu}$. The multiplication law in $G^g_M$ is
determined with the help of functions $\tilde \varphi^\alpha$ and
$\tilde f^\mu$ which determine the multiplication law in the Lie
group $G_M$ and its action on the manifold $M$ by the formulae
\cite{3}:
\begin{equation}
\label{eq2} (\tilde{g}\times\tilde{g}')^{\alpha}(x)
=\tilde{\varphi}^{\alpha}(\tilde{g}(x), \tilde{g}'(x')),
\end{equation}
\begin{equation}\label{eq3}
x'^{\mu}=\tilde{f}^{\mu}(x, \tilde{g}(x)).
\end{equation}

The formula \eqref{eq3} sets the action of $G^g_M$ on $M$.
\begin{definition}
The groups $G^g_M$, which are parameterized by smooth functions
$\tilde g^\alpha (x)$ with property \eqref{eq1} and have
multiplication law \eqref{eq2}, \eqref{eq3} are called the
\textbf{generalized gauge groups}.
\end{definition}

Let's pass from the group $G^g_M = \{\tilde g(x)\}$ to the group
$G^{gH}_M = \{g(x)\}$ isomorphic to it in accordance with the
formula
\begin{equation}\label{eq4}
g^a (x)=H^a_x (\tilde{g}(x))
\end{equation}
(Latin indices assume the same values as the corresponding Greek
ones). The smooth maps $H_x : G_M \rightarrow G_M$ have the
properties:
\[
1H)\ H_{x}(0)=0 \quad{} \forall x \in M ;
\]
\[
2H)\ \exists H^{-1}_x (g): \quad {} H^{-1}_x (H_x (g))= g \quad {}
\forall g \in G_M,\ x \in M.
\]

The multiplication law of the group $G^{gH}_M$ is determined by
its isomorphism \eqref{eq4} to the group $G^{g}_M$ and the
formulae \eqref{eq2} and \eqref{eq3}:
\begin{equation}
\label{eq5} (g \ast g')^a (x)=\varphi^a (x,g(x),g'(x')) :=
H^a_x(\tilde\varphi(H^{-1}_x(g(x)),H^{-1}_{x'}(g'(x')))),
\end{equation}
\begin{equation}
\label{eq6} x'^{\mu}=f^{\mu}(x,g(x)) := \tilde f^\mu
(x,H^{-1}_x(g(x))).
\end{equation}
The formula \eqref{eq6} sets the action of $G^{gH}_M$ on $M$.
\begin{definition}
Transformations \eqref{eq4} between the groups $G^g_M$ and
$G^{gH}_M$ are called the \textbf{deformations} of generalized
gauge groups and the groups $G^{gH}_M$ are called the
\textbf{infinite (generalized gauge) deformed groups} \cite{3}.
\end{definition}

In the set $D=\{H_x\}$ of maps $H_x$ the multiplication law can be
defined:
\begin{equation}\label{eq7}
(H_1 \circ H_2 )_x (g) :=H_{1x}(H_{2x}(g)) .
\end{equation}

The set $D$ becomes a group according to this multiplication law.
\begin{definition}
The maps $H_x : G_M \rightarrow G_M$ with properties 1$H$, 2$H$
are called the \textbf{deformation maps} (functions $H^a_x (g)$
are called the \textbf{deformation functions}), the group
$D=\{H_x\}$ of deformation maps with the multiplication law
\eqref{eq7} - the \textbf{group of deformations} \cite{3}.
\end{definition}

The functions $h(x)^a_\alpha := \partial_{\tilde{\alpha}} H^a_x
(\tilde g)|_{\tilde g=0}$, where $\partial_{\tilde{\alpha}}
:=\partial/\partial\tilde g^\alpha$, are called \emph{deformation
coefficients.}

With the help of the coefficients of the expansion
\begin{equation}
\label{eq8} \varphi^a (x,g,g')=g^a + g'^a + \gamma(x)^a{}_{bc}\
g^b g'^c + \frac{1}{2} {\rho(x)^a}_{bcd}\ g^d g'^b g'^c + \ldots
\end{equation}
functions
\begin{equation}\label{eq9}
F(x)^a_{bc} := \gamma(x)^a{}_{bc} - \gamma(x)^a{}_{cb} ,
\end{equation}
\begin{equation}
\label{eq10} R(x)^a{}_{dbc} := \rho(x)^a{}_{dbc} -
\rho(x)^a{}_{dcb}
\end{equation}
are defined, which are their skew-symmetric parts. They are called
the \emph{structure functions} (versus the structure constants for
ordinary Lie groups) and the \emph{curvature coefficients} of the
deformed group $G^{gH}_M$ respectively.

Since $\xi(x)_a^\mu := \partial_a f^\mu (x,g) |_{g = 0} =
h(x)^\alpha_a \tilde \xi (x)^\mu_\alpha$, where $\partial_b :=
\partial / \partial g^b$ and $h(x)^\alpha_a$ is
reciprocal to the $h(x)_\alpha^a$ matrix, the generators $X_a :=
\xi(x)_a^\mu
\partial_\mu$ ($\partial_\mu := \partial / \partial x^\mu$) of the
deformed group $G^{gH}_M$ are expressed with the help of
generators $\tilde X_a := \tilde\xi(x)_a^\mu
\partial_\mu$ of
the group $G_M$ and deformation coefficients: $X_a = h(x)_a^\alpha
\tilde X_\alpha$. So in an infinitesimal (algebraic) level,
deformation is reduced to nondegenerate linear transformations of
generators of the initial Lie group independent in every point $x
\in M$.
\begin{proposition}
Generators commutators of the deformed group $G^{gH}_M$ are the
linear combinations of generators with structure functions which
are the coefficients \cite{3}:
\begin{equation}\label{eq11}
[X_a , X_b ] = F(x)^c_{ab} X_c .
\end{equation}
\end{proposition}

The equation \eqref{eq11} generalize the Maurer-Cartan equation
\begin{equation}\label{eq12}
[\tilde X_\alpha , \tilde X_\beta ] = \tilde
F^\gamma_{\alpha\beta} \tilde X_\gamma
\end{equation}
for the infinite deformed groups $G^{gH}_M$, where $\tilde
F^\gamma_{\alpha\beta}$ are the structure constants of the initial
Lie group $G_M$. The equation \eqref{eq11} is reduced to the
equation \eqref{eq12} for the generalized gauge nondeformed group
$G^g_M$.

%%%%%%%%%%%%%%%%%%%%%%%%%%%%%%%%%%%%%%%%%%%%%%%%%%%%%%%%%%%%%%%%
\section{The Deformed Group of Diffeomorphisms \\
and Geometrical Structure of Riemannian Space} \protect\label{DGD}
%%%%%%%%%%%%%%%%%%%%%%%%%%%%%%%%%%%%%%%%%%%%%%%%%%%%%%%%%%%%%%%%

Let $G_M = T_M$, where $T_M$ is the group of translations. In this
case $(\tilde t\cdot\tilde t')^{\mu}=\tilde t^{\mu}+\tilde t'^\mu
$ and $x'^{\mu}=x^{\mu}+ \tilde t^\mu$. The group $T^g_M$ is
parameterized by the functions $\tilde t^{\mu}(x)$, which satisfy
the condition $\det\{\delta_{\nu}^{\mu}+\partial_{\nu}\tilde
t^{\mu}(x)\}\neq0$, $\forall x\in M$. The multiplication law in
$T^g_M$ is
\begin{equation}
\label{eq13} (\tilde t \times \tilde t')^\mu (x) = \tilde t^\mu
(x) + \tilde t'^\mu (x'),
\end{equation}
\begin{equation}
\label{eq14} x'^\mu = x^\mu + \tilde t^\mu (x) ,
\end{equation}
where \eqref{eq14} determines the action of $T^g_M$ on $M$. The
multiplication law indicates that $T^g_M$ is the group of
diffeomorphisms $Diff\ M$ in additive parametrization. The
generators of the  $T^g_M$-action \eqref{eq14} on $M$ are simply
derivatives $\tilde X_\mu =\partial_\mu$ and this fact corresponds
to the case of the flat space $M$.

Suppose that the group $T_{M}^{g}$ is deformed
$T_{M}^{g}\rightarrow T_{M}^{gH}$: $t^m(x)=H^m_x (\tilde{t}(x))$.
The multiplication law in $T^{gH}_M$ is determined by the
formulae:
\begin{equation}
\label{eq15} (t \ast t')^m (x)=\varphi^m (x,t(x),t'(x')) := H^m_x
(H^{-1}_x (t(x))+H^{-1}_{x'} (t'(x'))),
\end{equation}
\begin{equation}
\label{eq16} x'^{\mu}=f^{\mu}(x,t(x)) := x^\mu + H^{-1\mu}_x
(t(x)).
\end{equation}
Formula \eqref{eq16} sets the action of $T^{gH}_M$ on $M$.

Let's consider expansion
\begin{equation}
\label{eq17} H^{m}_x (\tilde{t})=h(x)_{\mu}^{m}[\tilde{t}^{\mu}+
\frac{1}{2}\Gamma(x)_{\nu\rho}^{\mu}\tilde{t}^{\nu}\tilde{t}^{\rho}+
\frac{1}{6}\Delta(x)_{\nu\rho\sigma}^{\mu}\tilde{t}^{\nu}\tilde{t}^{\rho}
\tilde{t}^{\sigma}+\ldots].
\end{equation}
Using of the formula \eqref{eq15}, for coefficients of expansion
\eqref{eq8} we can obtain
\begin{equation}
\label{eq18} \gamma^{m}{}_{kn}=h_{\mu}^{m}(\Gamma_{kn}^{\mu}+
h_{k}^{\nu}\partial_{\nu}h_{n}^{\mu}),
\end{equation}
\begin{equation}
\label{eq19} \rho^{m}{}_{lkn}=h_{\mu}^{m}(\Delta_{lkn}^{\mu}-
\Gamma_{ns}^{\mu}\Gamma_{kl}^{s}-
h_{n}^{\nu}\partial_{\nu}\Gamma_{\kappa\lambda}^{\mu}
h_{k}^{\kappa}h_{l}^{\lambda}).
\end{equation}
So formulae \eqref{eq9} and \eqref{eq10} for the structure
functions and the curvature coefficients of deformed group
$T^{gH}_M$ yield
\begin{equation}
\label{eq20}
F_{\mu\nu}^{n}=-(\partial_{\mu}h_{\nu}^{n}-\partial_{\nu}h_{\mu}^{n}),
\end{equation}
\begin{equation}
\label{eq21} R^{\mu}{}_{\lambda\kappa\nu}=
\partial_{\kappa}\Gamma_{\nu\lambda}^{\mu}-\partial_{\nu}\Gamma_{\kappa\lambda}^{\mu}
+\Gamma_{\kappa\sigma}^{\mu}\Gamma_{\nu\lambda}^{\sigma}-
\Gamma_{\nu\sigma}^{\mu}\Gamma_{\kappa\lambda}^{\sigma}.
\end{equation}
In this formulae matrix $h^m_\mu$ and reciprocal to it matrix
$h_m^\mu$ we use for changing Greek indices to Latin (and vice
versa).

Formulae \eqref{eq20} and \eqref{eq21} show that groups $T^{gH}_M$
contain the information about the geometrical structure of the
space $M$ where they act. The generators $X_k = h^\nu_k
\partial_\nu$ of the $T^{gH}_M$-action \eqref{eq16} on $M$ can be
treated as \emph{affine frames}. Structure functions
$F^n_{\mu\nu}$ differ from the \emph{anholonomity coefficients}
only by the factor $- 1/2$.

Let us write the multiplication law of the group $T^{gH}_M$
\eqref{eq15} for the infinitesimal second factor:
\begin{equation}
\label{eq22} (t \ast \tau)^m (x)=t^m
(x)+\lambda(x,t(x))^{m}{}_{n}\ \tau^{n}(x'),
\end{equation}
where
$\lambda(x,t)^{m}{}_{n}:=\partial_{n'}\varphi^{m}(x,t,t')\mid_{t'=0}$.
Formula \eqref{eq22} gives the rule of the addition of vectors,
which are set in different points $x$ and $x'$ or \emph{the rule
of the parallel transport} of the vector field $\tau$ from point
$x'$ to point $x$:
\begin{equation}
\label{eq23} \tau_{\parallel}^{m}(x) =\lambda(x,t(x))^{m}{}_{n}\
\tau^{n}(x')
\end{equation}
or in the coordinate basis
\begin{equation}\label{eq24}
\tau_{\parallel}^{\mu}(x) =\partial_{\tilde{\nu}}H^{\mu}_x (
\tilde{t}) \tau^{\nu}(x+\tilde{t})
\end{equation}
where $\partial_{\tilde{\nu}} :=\partial/\partial\tilde{t}^\nu$.
This formula determines the \emph{covariant derivative}
\begin{equation}
\label{eq25}
\nabla_{\nu}\tau^{\mu}(x)=\partial_{\nu}\tau^{\mu}(x)+
\Gamma(x)_{\sigma\nu}^{\mu}\tau^{\sigma}(x),
\end{equation}
where functions $\Gamma(x)^\mu_{\sigma\nu}$ seting the second
order of the expansion \eqref{eq17} of deformation functions, play
the role of \emph{coefficients of an affine connection} in the
coordinate basis. So, curvature coefficients \eqref{eq21}
$R^{\mu}{}_{\lambda\kappa\nu}$ of the group $T^{gH}_M$ coincide
with the \emph{Riemann curvature tensor.} The functions
$\Gamma(x)^\mu_{\sigma\nu}$ are symmetric on the bottom indices,
so torsion is equal to zero. Relationship \eqref{eq18} means that
coefficients $\gamma^m{}_{kn}$, which set the second order of the
expansion \eqref{eq8} of the multiplication law in the group
$T^{gH}_M$, are coefficients of the \emph{affine connection in the
affine basis} $X_k$.

At the consecutive performance of the deformations $H_{2x}$ and
$H_{1x}$ for the resulting deformation $H_{3x}=(H_1 \circ H_2 )_x$
one can obtain:
\begin{equation}\label{eq26}
h_3{}^m_\mu =h_1{}^m_p\ h_2{}^p_\mu ,\quad\Gamma_3{}^m_{\mu\nu}
=\Gamma_1{}^m_{ps}\ h_2{}^p_\mu\ h_2{}^s_\nu +h_1{}^m_p\
\Gamma_2{}^p_{\mu\nu}.
\end{equation}
The last formula corresponds to the notion of deformations of
connections \cite{4} and explains the term "deformations" in our
case.

Suppose now, that generators $X_k =h^\nu_k \partial_\nu$ of the
$T^{gH}_M$-action on $M$ \eqref{eq16} are orthonormalized frames,
i.e. $g(X_m,X_n)=\eta_{mn}$, where $\eta_{mn}$ - Euclidean metric,
and infinitesimal parallel transports of vector fields lead only
to their rotations, i.e.:
\begin{equation}\label{eq27}
\lambda(x,t)^{m}{}_n \cong \delta^m_n + \gamma^{m}{}_{kn}\ t^k \in
SO(n).
\end{equation}
For coefficients $\gamma^{m}{}_{kn}$ this gives
\begin{equation}\label{eq28}
\gamma^\cdot_{ksl} +\gamma^\cdot_{lsk} =0,
\end{equation}
(we fulfill lowering indices with the help of the metric:
$\gamma^\cdot_{mkl}:=\eta_{mn}\gamma^{n}{}_{kl}$). Together with
formula \eqref{eq18} the equation \eqref{eq28} gives the condition
of coordination of connection with the metric
$g(\partial_\mu,\partial_\nu)=:g_{\mu\nu} =h^m_\mu h^n_\nu
\eta_{mn}$:
\begin{equation}\label{eq29}
\Gamma^\cdot_{\mu\nu\sigma} +\Gamma^\cdot_{\nu\mu\sigma}
=\partial_\sigma g_{\mu\nu} ,
\end{equation}
With the condition of torsion vanishing, this yields that
coefficients $\Gamma^\rho_{\mu\nu}$ may be written as
\begin{equation}
\label{eq30}
\Gamma_{\mu\nu}^{\rho}=\frac{1}{2}g^{\rho\sigma}(\partial_{\mu}g_{\nu\sigma}
+\partial_{\nu}g_{\mu\sigma}-\partial_{\sigma}g_{\mu\nu}).
\end{equation}
So these coefficients coincide with the \emph{Christoffel symbols}
$\{{}^\rho_{\mu\nu}{} \}$.
\begin{proposition}
The deformed group $T^{gH}_M$ of diffeomorphisms of coordinate
chart $U \subset M$, which is obtained with satisfying the
condition \eqref{eq28} (or \eqref{eq29}), acting on $U$ sets on it
structure of a Riemannian space. Geometrical characteristics of
the space $U$ (connection coefficients, curvature tensor etc.) are
contained in the multiplication law of the group $T^{gH}_M$. Thus
any Riemannian structure on $U \subset M$ may be set \cite{2}.
\end{proposition}

This proposition realizes Klein's Erlangen Program for the
Riemannian space.

%%%%%%%%%%%%%%%%%%%%%%%%%%%%%%%%%%%%%%%%%%%%%%%%%%%%%%%%%%%%%%%%
\section{The Parallel Transports as Automorphisms \\
of Deformed Groups of Diffeomorphisms}
\protect\label{ADGD}
%%%%%%%%%%%%%%%%%%%%%%%%%%%%%%%%%%%%%%%%%%%%%%%%%%%%%%%%%%%%%%%%

In the approach considered in previous section the vector fields
in the curved Riemannian space were presented by infinitesimal
parameters of the deformed group of diffeomorphisms. So, we can
consider the parallel transports of vector fields \eqref{eq23} as
certain automorphisms of the deformed group of diffeomorphisms
$T^{gH}_M$.

Let's consider a transformation
\begin{equation}
\label{eq31} \tau_{\parallel}^{m}(x) :=\tilde H^m_x ((t\ast
\tau)(x)-t(x))=\tilde H^m_x (\varphi(x,t(x),\tau(x'))-t(x))
\end{equation}
where
\begin{equation}
\label{eq32} x'^\mu =f^\mu(x,t(x)):=x^\mu +H^{-1\mu}_x (t(x)),
\end{equation}
$\tilde H_x$ are variable deformation maps, and $H_x$ - the fixed
deformation map, which was used for constructing the group
$T^{gH}_M$ and which defines geometrical characteristics of space
where we intend to consider parallel transports. The inverse for
\eqref{eq31} transformation is:
\begin{equation}
\label{eq33} \tau^m (x)=\varphi^m (x,t^{-1}(x),t(\tilde x)+\tilde
H^{-1}_{\tilde x} (\tau_\| (\tilde x))) ,
\end{equation}
where $\tilde x^\mu :=f^\mu (x,t^{-1}(x))$.

Let's pass from the group $T^{gH}_M =\{\tau(x)\}$ to the group
$T_\|{}^{gH}_M =\{\tau_\| (x)\}$, isomorphic to it by the formula
\eqref{eq31}.

The multiplication law in the group $T_\|{}^{gH}_M$ is determined
by its isomorphism \eqref{eq31} to the group $T^{gH}_M$ and by the
multiplication law \eqref{eq15}, \eqref{eq16} in the group
$T^{gH}_M$. The group $T_\|{}^{gH}_M$ acts on the chart $U$
according to the formula
\begin{equation}
\label{eq34} \bar x^\mu =f^\mu_\| (x,\tau_\| (\tilde x)):= f^\mu
(x,\tau (x)) = f^\mu (\tilde x,t(\tilde x) + \tilde
H^{-1}_{\tilde{x}} (\tau_\|(\tilde x))) .
\end{equation}
We should emphasize that the transformation of a point $x$ is
determined by the value of functions $\tau_\| (\tilde x)$ (that
parameterizes the group of parallel transports $T_\|{}^{gH}_M$) in
another point $\tilde x$.

We shall name the group $T^{gH}_M$ with infinitesimal parameter
$\tau(x)$ an \emph{infinitesimal group} $T^{gH}_M$. For
infinitesimal group $T^{gH}_M$ from \eqref{eq31} and \eqref{eq34}
follows:
\begin{equation}
\label{eq35} \tau_{\parallel}^{m}(x) =L(x)^m_p
\lambda(x,t(x))^{p}{}_{n}\ \tau^{n}(x') ,
\end{equation}
\begin{equation}
\label{eq36} h_{\parallel}(x)_{m}^\mu =h(x)^\mu_k \lambda^{-1}
(\tilde x,t(\tilde x))^{k}{}_{n}\ L^{-1} (\tilde x)^{n}_m ,
\end{equation}
where in this case $L(x)^m_p :=\partial_p \tilde H^m_x (t)|_{t=0}$
and $h_\parallel (x)^\mu_m := \partial_m f^\mu_\parallel
(x,\tau)|_{\tau =0} .$

Transformations \eqref{eq35} (or \eqref{eq36}) form an infinite
group $DT$ with parameters $g(x)=\{t(x),L(x)\}$ and multiplication
law
\begin{equation}
\label{eq37} (g\ast g')^m (x) =\varphi^m (x,t(x),t'(x')),
\end{equation}
\begin{equation}
\label{eq38} (g\ast g')^m_n (x)=L(x)^m_p
\lambda(x,t(x))^{p}{}_{r}\ L'(x')^r_s
\lambda(x',t'(x'))^{s}{}_{t}\ \lambda^{-1} (x,\varphi
(x,t(x),t'(x')))^{t}{}_n ,
\end{equation}
where
\begin{equation}
\label{eq39} x'^\mu =f^\mu (x,t(x))
\end{equation}
and we consider that $g^m (x)=t^m (x)$ and $g^n_m (x)=L(x)^n_m$.
This multiplication law shows, that group $DT$ has the structure
$T^{gH}_M \times ) GL^g (n)$, where $T^{gH}_M
=\{t(x),\lambda^{-1}(x,t(x))\}$ and $GL^g (n)=\{0,L(x)\}$ are its
subgroups. Moreover, the group $DT$ is the deformed generalized
gauge group $(T\otimes GL(n))^{gH}_M$ \cite{3}.

Formula \eqref{eq39} determines the action of the group $DT$ on
the chart $U \subset M$, formula \eqref{eq35} on tangent vectors
and \eqref{eq36} on affine frames $X_m =h^\mu_m
\partial_\mu$ over $U$ respectively.
\begin{definition}
The group $DT=\{t(x),L(x)\}$ of automorphisms \eqref{eq35} of the
infinitesimal deformed group of diffeomorphisms $T^{gH}_M$ with
the multiplication law \eqref{eq37}-\eqref{eq39}, which act on
tangent vectors and affine frames over $U \subset M$ according to
formulae \eqref{eq35} and \eqref{eq36} respectively is called the
\textbf{group of parallel transports} in the space $U$.
\end{definition}

Let's consider structure functions $F(x)^c_{ab}$ of the group of
parallel transports $DT$ on the condition if $a=k$, $b=l$ (that
corresponds to the translation parameters). For $c=m$ from formula
\eqref{eq37} we obtain
\begin{equation}
\label{eq40} F^m_{kl} =h^m_\mu (h^\nu_k \partial_\nu h^\mu_l
-h^\nu_l \partial_\nu h^\mu_k )
\end{equation}
and for $c={}^{m}_{n}$ from formula \eqref{eq38} -
\begin{equation}
\label{eq41} F^{m}_{n}{}_{kl} =-\gamma^{m}{}_{sn}\
F^s_{kl}+h^\sigma_k
\partial_\sigma \gamma^{m}{}_{ln}-h^\sigma_l
\partial_\sigma \gamma^{m}{}_{kn}+\gamma^{m}{}_{ks}\
\gamma^{s}{}_{ln}-\gamma^{m}{}_{ls}\ \gamma^{s}{}_{kn} .
\end{equation}
These equations show that the structure functions $F^m_{kl}$ and
$F^m_{n}{}_{kl}$ of the group of parallel transports $DT$ coincide
with the structure functions $F^m_{kl}$ and curvature coefficients
$R^m {}_{nkl}$ of the deformed group of diffeomorphsms $T^{gH}_M$,
i.e. with anholonomity coefficients (with the factor $-2$) and the
Riemann curvature tensor (written in the affine frame)
respectively.

Generators $X^\tau_a$ of the action \eqref{eq35} of the group $DT$
on the tangent vectors for $a={}^{m}_{n}$ are
$(X^\tau{}^{n}_{m})^k_l =\delta^k_m \delta^n_l$ and for $a=m$ are
\begin{equation}
\label{eq42}(X^\tau_{m})^k_l =X_m \delta^k_l +\gamma^{k}{}_{ml}
\end{equation}
and coincide with covariant derivatives $X^\tau_m =\nabla_m$ in
the affine frame.

From the the generalized Maurer-Cartan equation \eqref{eq11} for
the group of parallel transports $DT$ follows the equation
\begin{equation}
\label{eq43}[\nabla_k ,\nabla_l ]^m_n =F^s_{kl} \nabla^m_{s}{}_{n}
+R^{m}{}_{nkl} ,
\end{equation}
which is equivalent to the structure equations of the curved space
of the torsion-free affine connection with the variable curvature
$R^{m}{}_n$ (if this connection satisfies condition \eqref{eq28}
of the Riemannian space):
\begin{equation}
\label{eq44} d\omega^m =\omega^n \wedge \omega^{m}{}_n ,
\end{equation}
\begin{equation}
\label{eq45} d\omega^{m}{}_n =\omega^{k}{}_n \wedge \omega^{m}{}_k
+R^{m}{}_{n} ,
\end{equation}
where $\omega^m =h^m_\mu dx^\mu$, $\omega^{m}{}_n
=\gamma^{m}{}_{\mu n}\ dx^\mu$ and $R^{m}{}_n
=\frac{1}{2}R^{m}{}_{n\mu\nu}\ dx^\mu \wedge dx^\nu$.

Formula \eqref{eq43} indicates that $R^m{}_{nkl} =0$ is the
necessary and sufficient condition that the set of translations
$\{t(x),1\}$ in the group $DT$ should form a subgroup. Gauge
linear transformations $L(x)$ in the case of a curved Riemannian
space are necessary for ensuring the group structure of the group
$DT$.

Formula \eqref{eq36} describes motion of the mobile frame
$X_{\parallel_m} =h^\mu_{\parallel_m}\partial_\mu$ at the
transformations of parallel transports from the group $DT$. For
infinitesimal translations (and finite linear transformation
$L^m_n$) the formula \eqref{eq36} gives
\begin{equation}
\label{eq46} X_{\parallel_m} =\bar{X}_m -\bar{t}^s
\bar{\gamma}^{n}{}_{sm}\ \bar{X}_n ,
\end{equation}
where $\bar{X}_m =L^{-1n}{}_m\ X_n$, $\bar{t}^s=L^{s}{}_n\ t^n$
and
\begin{equation}
\label{eq47} \bar{\gamma}^{n}{}_{sm}=L^{n}{}_l (\gamma^{l}{}_{rn}\
L^{-1r}{}_{s}\ L^{-1n}{}_{m}+L^{-1n}{}_{s}\ h^\sigma_n
\partial_\sigma L^{-1l}{}_{m}).
\end{equation}
Formula \eqref{eq46} can be used for the definition of covariant
derivatives in the mobile frame terms:
\begin{equation}
\label{eq48} \nabla_{X_s}X_{\parallel_m} =\lim_{t^{s}\rightarrow
0} (X_m -X_{\parallel_m})/t^s =\gamma^{n}{}_{sm}\ X_n .
\end{equation}

Let's suppose that the group $T^{gH}_M$ is obtained with the
fulfilling condition \eqref{eq28} (or \eqref{eq29}). In this case
we can prove the next proposition.
\begin{proposition}
Parallel transports of vector fields in curved Riemannian space
are described by the group $DT$ of special automorphisms of the
infinitesimal deformed group of diffeomorphisms $T^{gH}_M$.
Translations generators of the group $DT$ are the covariant
derivatives of vector fields, and the structure functions of the
group $DT$ contain the curvature tensor.

The equations of structure of Riemannian space \eqref{eq44},
\eqref{eq45} (Cartan equations) are the necessary and sufficient
conditions of the group $DT$ existence.
\end{proposition}

The group $DT$, as well as the group $T^{gH}_M$, contains
information about the structure of the Riemannian space on
$U\subset M$. Generators and structure functions of group $DT$
contain this information. The structure of Riemannian space is set
on $U$ at infinitesimal action of group $DT$ in the tangent bundle
of space $U$ while for setting of the Riemannian structure on $U$
with the help of the group $T^{gH}_M$ it is necessary to consider
its action on $U$, at least, up to the second order on
translations $t$ inclusively.

%%%%%%%%%%%%%%%%%%%%%%%%%%%%%%%%%%%%%%%%%%%%%%%%%%%%%%%%%%%%%%%%
\subsection*{Acknowledgments}
%%%%%%%%%%%%%%%%%%%%%%%%%%%%%%%%%%%%%%%%%%%%%%%%%%%%%%%%%%%%%%%%

We wish to thank professor M.I. Jaloviy for his assistance in
preparing the English version of the manuscript.

%%%%%%%%%%%%%%%%%%%%%%%%%%%%%%%%%%%%%%%%%%%%%%%%%%%%%%%%%%%%%%%%

%%%%%%%%%%%%%%%%%%%%%%%%%%%%%%%%%%%%%%%%%%%%%%%%%%%%%%%%%%%%%%%%

\end{document}